\def\**{$*\!*\!*$}


 at 12truept
\font\tenmsy=msym10
\textfont8=\tenmsy
\mathchardef\ssm="7872

\input psfig
\mathsurround = 2pt
\abovedisplayskip=6pt
\belowdisplayskip=6pt

\def \QED {\rlap{$\sqcup$}$\sqcap$\smallskip}

\def\ref{\hangindent=1pc \hangafter=1 \noindent}

\input psfig

\centerline{\bf Teichm\"{u}ller  space of
           Fibonacci maps.}\smallskip

\centerline{March 18, 1993}
\medskip

\centerline{Mikhail Lyubich}\smallskip
\centerline{Mathematics Department and IMS, SUNY Stony Brook}\bigskip

\noindent {\bf \S1. Introduction.}
According to Sullivan, a space ${\cal E}$
of unimodal maps with the same combinatorics
(modulo smooth conjugacy) 
should be treated as an infinitely-dimensional Teichm\"{u}ller space.
This is a basic idea in Sullivan's approach to the Renormalization Conjecture
[S1], [S2].  One of its principle ingredients is to supply ${\cal E}$ 
with the Teichm\"{u}ller metric.
To have such a metric one has to know, first of all, that all maps
of ${\cal E}$ are quasi-symmetrically conjugate.
This was proved in [Ji] and [JS]    for some classes of 
non-renormalizable maps (when the critical point is not too recurrent).
Here we consider a space of non-renormalizable unimodal  maps
 with in a sense fastest possible recurrence of the critical point
(called Fibonacci). 
 Our goal is to supply this space 
with the Teichm\"{u}ller metric.

Let $f$ be a unimodal map with critical point $c$. 
A Fibonacci unimodal map $f$ can be defined by saying that the closest
returns of the critical point occur at the Fibonacci moments.
This  combinatorial type was suggested by Hofbauer and Keller [HK]
as an extremal among non-renormalizable types (see [LM] for more 
detailed history).  Its combinatorial, geometric
and measure-theoretical properties were studied in
[LM]  under the assumptions that $f$ is {\sl quasi-quadratic}, i.e., it is
 $C^2$-smooth and has the quadratic-like
critical point (see also [KN]). We will assume this 
regularity throughout the paper.  

A principle object of our combinatorial considerations is a nested sequence
of intervals $I^0\supset I^1\supset...$ obtained subsequently by
pulling back along the critical orbit. 
Our proof is based upon the geometric result of [LM] which says that
the scaling factors $\mu_n=|I^n|/|I^{n-1}|$
characterizing the geometry of the Fibonacci map
decay exponentially. It follows that appropriately defined
renormalizations $R^n f$ are becoming purely quadratic near the critical point.
This reduces the renormalization process to the iterates of quadratic maps.

The next idea is to consider a quasi-conformal continuation of $f$ to the
complex plane which is asymptotically conformal on the real line. 
Then we consider complex generalized renormalizations, and prove that
the renormalized maps are becoming purely quadratic in the complex plane
as well.  Hence the geometric patterns of  renormalized
maps are subsequently obtained by the Thurston pull-back transformation
(up to an exponentially small error) in an appropriate 
 Teichm\"{u}ller space. It follows
that these patterns converge (after rescaling) to the corresponding pattern 
of the quadratic map $p: z\mapsto z^2-1$. In particular, the shape of the
complex puzzle-pieces converge to the Julia set of $p$, 
see Figure~1
(this is perhaps the most unexpected result of our analysis).   
  
\midinsert{
 \centerline{\psfig{figure=puzzle1.ps,width=.75\hsize}}
 \smallskip
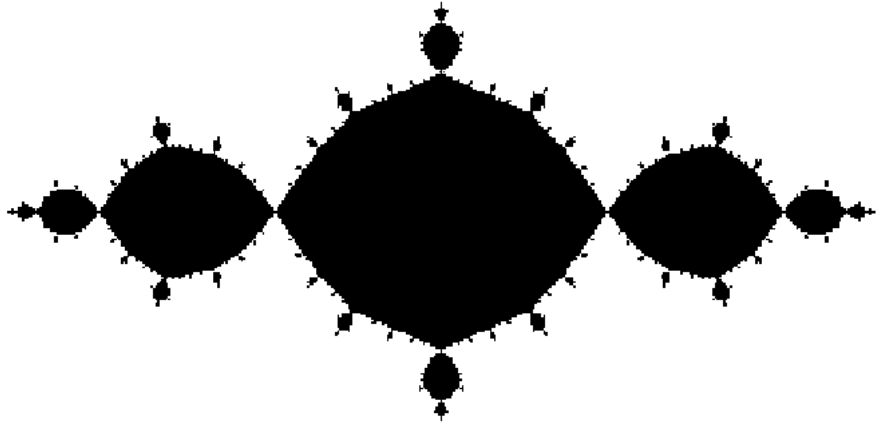
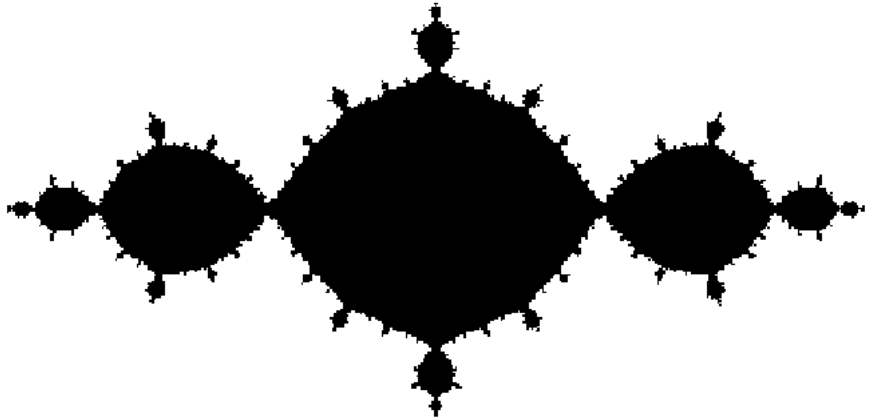
 \centerline{Figure 1. A Fibonacci puzzle-piece (below) versus the Julia set
            of $z\mapsto z^2-1$. }
 \centerline{(made by S.~Sutherland and B.~Yarrington)}
}
\endinsert


To each renormalization we then associate a pair of pants $Q_n$
by removing from the critical puzzle-piece of level $n$ two puzzle-pieces
of the next level. 
Using a same type of argument as above,
we show that the pairs of pants $Q^n$ and $\tilde Q^n$ stay on 
bounded distance.
This yields the quasi-conformal equivalence of the critical sets of $f$
and $\tilde f$. 

To complete the construction of the quasi-symmetric conjugacy,
 we apply  a Sullivan-like pull-back argument. 
However, this is not quite straightforward
since there is no 
dilatation control  away from the real line.  
  
In the last section we prove that two Fibonacci maps which stay on zero
Teichm\"{u}ller distance are smoothly conjugate. So this pseudo-metric
is non-degenerate on the smooth equivalence classes. 

We will use abbreviations qc and qs for ``quasi-conformal" and 
``quasi-symmetric" respectively.

\medskip\noindent 
{\bf Remark 1.} Since the rate at which the scaling factors decrease depends
on the initial bounds of the map only, the dilatation of the conjugacy
we construct also depends only on this data.

\medskip\noindent 
{\bf Remark 2.} 
It is proved in [L] that, as in the Fibonacci case, 
the scaling factors of any non-renormalizable quasi-quadratic map
 decay exponentially.  This allows us to generalize the above result to
all combinatorial classes of quasi-quadratic maps. 
The exposition of this result
 is more technical, and it will be the subject of forthcoming notes.
Note that for polynomial-like maps this result follows from the Yoccoz
Theorem (see [H] for the exposition of this theorem, and [K] for an
alternative proof  based upon a pull-back argument).

\medskip\noindent 
{\bf Remark 3.} In this paper we concentrate on the dynamical constructions,
and don't touch the issue of the sharp regularity
for which the theory can be built up. This issue is clearly important for
a proper Teichm\"{u}ller theory (compare [S2] and [G]), 
and will be discussed elsewhere. 
\bigskip

\noindent{\bf \S2. Asymptotically conformal continuation and 
generalized  renormalization.}
  
\noindent{\bf Real renormalization} (see [LM]). Given a Fibonacci map $f$,
there is  a sequence of maps 
$$g_n: I^n_0\cup I^n_1\rightarrow I^{n-1}_0,\quad n=1,2,...$$
constructed in the following way. Let $I^0\equiv I^0_0$ be a $c$-symmetric
interval satisfying the property $f^n(\partial I^0)\cap I^0=\emptyset$,
$n=1,2,...$.
Now given $I^{n-1}\equiv I^{n-1}_0\ni c$ by induction, 
let us consider the first return map
$f_n: \cup I^n_j\rightarrow I^{n-1}$.
 Its domain of definition generally
consists of infinitely many intervals $I^n_j\subset I^{n-1}$.
 However, for the Fibonacci map
only two of them, $I^n\equiv I^n_0\ni c$ (the ``central" one)
and $I^n_1$ intersect the critical set $\omega(c)$.
Let us define $g_n$ as the restriction of $f_n$ to these two intervals.
 These maps satisfy the following properties:
\item (i) $g_n: I^n_1\rightarrow I^{n-1}_0$ is a diffeomorphism and
 $g_n (\partial I^n_0)\subset \partial I^{n-1}_0$; 
\item (ii)  $g_n I^n_0\supset I^n_0$ ({\sl high return});
\item (iii) $g_n c\in I^n_1$ and $g^2_n c\in I^n_0$.

By rescaling $I^n$ to some definite size $T$ (e.g., $T=[0,1]$), we
obtain the generalized $n$-fold renormalization
$$R^n f: T^n_0\cup T^n_1\rightarrow T$$ of $f$.
The asymptotic properties of the renormalized maps express the small
scale information of the critical set $\omega(c)$.

Let us now introduce the principle geometric parameters, the scaling factors
$$\mu_n={|I^n|\over |I^{n-1}|}={|T^n|\over |T|}.$$
The main result of [LM] says that they decrease to 0 exponentially at the
following rate:
$$\mu_n\sim a\left({1\over 2}\right)^{n/3}.  \eqno (1)$$
It follows by the Koebe principle that up to an exponentially small error
the restriction of $R^n f$ to the central interval $T^n_0$ is purely quadratic,
 while the restriction to $T^n_1$ is linear. This all we need to know for
the comprehensive study of $f$.

\medskip\noindent
{\bf Asymptotically conformal continuation.} Let us represent $f$ as
$h\circ \phi$ where $\phi(z)=(z-c)^2$ is the quadratic map, while
 $h$ is a $C^2$-diffeomorphism of appropriate intervals.
Let us continue $h$ to a diffeomorphism of a bounded $C^2$ norm on
 the whole real line, and then
 consider the Ahlfors-Beurling continuation of $h$ to the complex plane:
$$\hat{h}(x+iy)={1\over 2y}\int_{x-y}^{x+y} h(t) dt+
      {1\over y}\left(\int_x^{x+y} h(t) dt - \int_{x-y}^x h(t) dt\right).$$
This is clearly a $C^2$-map, and one can check by calculation that
$\bar\partial \hat h=0$ on the real line.  Hence
$\bar\partial \hat h/\partial \hat h=O(|y|)$ as $|y|\to 0$.
This provides us with a $C^2$
extension of $f$ which is asymptotically conformal on the real line in 
the sense that
 $$\mu(z)\equiv\bar\partial \hat f/\partial \hat f=O(|y|) \eqno (2)$$
as well. In what follows we denote the extended $h$ and $f$ by the same
letters.

\medskip\noindent
{\bf Complex pull-back.} Given an interval $I\subset {\bf R}$ and 
$\theta\in (0,\pi/2)$, let 
$D_{\theta}(I)$ denote the domain bounded by the
 union of two ${\bf R}$-symmetric arcs of the circles
which touch the real line at angle $\theta$. In particular, 
$D_{\pi/2}(I)\equiv D(I)$ is the
Euclidean disk with diameter $I$.  Observe that $I$ is a hyperbolic geodesic
in the domain ${\bf C}\ssm({\bf R}\ssm I)$ and $D_{\theta}(I)$ is its
hyperbolic neighborhoods of radius depending only on $\theta$.

We say than an interval $\tilde I$ is
obtained from the $I$ by $\alpha$-scaling  if these intervals are
cocentric and $|\tilde I|=(1+\alpha)|I|$.

\proclaim Lemma 1. Let $\alpha<1$, $n$ be sufficiently big.
 Let us consider the $\alpha$-scaled interval 
$\tilde I^n\supset I^n$. Let $\Delta=D(\tilde I^n)$, and $\Delta'$ be the
pull-back of $\Delta$ by $g_{n+1}|I^{n+1}$. Then 
$\Delta'\subset D(\tilde I^{n+1})$ where $\tilde I^{n+1}$ is obtained from
$I^{n+1}$ by $\beta$-scaling with $\beta=\alpha+O(\mu_n)$.

\noindent{\bf Proof.} Let us skip the index $n$ in the notations of objects of
level $n$, while mark the objects of level $n+1$ with prime. Set $g|I'=f^p$,
and let us consider the pull back $I, I_{-1},...,I_{-p}\equiv I'$ of $I$ along 
the orbit $\{f^k c\}_{k=0}^p$.  Then 
$$\sum_{k=0}^p | I_{-k}|=O(\mu).   \eqno (3)$$
Since the map $f^k: I_{-k}\rightarrow I$ has the Koebe space covering
$I^{n-1}$, the  pull-back $\tilde I_{-k}$ 
of $\tilde I$ along the same orbit
also has the total length $O(\mu)$. 
 
Let us now take the disk $\Delta$ and pull it back along the same orbit.
We obtain a sequence of pieces $\Delta_{-k}$ based upon the intervals
$\tilde I_{-k}$. 
Assume by induction that $\Delta_{-l}\subset D_{\theta(k)}(\tilde I_{-l})$,
$l=0,...,k<p$, with
$$\theta_l=\alpha+O(\sum_{j=0}^{l-1} |\tilde I_{-j}|).  \eqno (4)$$

 Represent $f$ as $h\circ \phi$ and carry out the next pull back
in two steps: first by the diffeomorphism $h$ and then by the quadratic map 
$\phi$. Let $h^{-1} \tilde I_{-k}=L_{-k}$. If we rescale the intervals 
$\tilde I_{-k}$
and $L_{-k}$ to the unit size, the $C^1$-distance from
 the rescaled map $H^{-1}: [0,1]\rightarrow [0,1]$
  to id is $O(|I_{-k}|)$. It follows that
$$h^{-1} \Delta_{-k}\subset D_{\theta(k+1)}(L_{-k})  \eqno (5)$$
with $\theta(k+1)$ as in (4). 

Consider now two cases. Let first $k<p-1$. Then 
$\phi: L_{-k}\rightarrow\tilde I_{-(k+1)} $ is a diffeomorphism
and by the Schwarz lemma (see the above hyperbolic interpretation
of the $D_{\theta}(I)$)
$$\Delta_{-(k+1)}\subset D_{\theta(k+1)}(\tilde I_{-(k+1)}).$$

Let us now carry out the last pull-back corresponding to
$k=p-1$. Then $\phi|I_{-(k+1)}=\phi|I'$ is the quadratic folding map
into $L\equiv L_{-(p-1)}$.
Moreover, what is important is that $\phi I'$ covers at most half
(up to an error of order $O(\mu)$) of the interval $L$
(It follows from the high return property  of $g$ and the estimate of
its non-linearity). Hence we can find an interval $K\supset L$
centered at the critical value $gc$ such that
$$D_{\theta(p-1)}(L)\subset D(K)$$
and
$$|K|=2|\phi I'|(1+O(\mu)).$$
Two last equations together with (4) yield the required.  \QED

Let us now take the Euclidean disk $\Delta=D(I^m)$ and  pull it back
by  the maps $g_n$ continued to the complex plane. Denote the corresponding
domains by $\Delta^n_0$ and $\Delta^n_1$, $n>m$.

\proclaim Corollary 2. If $m$  is sufficiently big then 
the diam$\Delta^n_j$ is commensurable with the
diam$I^n_j$.

\noindent{\bf Proof.} Applying the previous lemma $n-m$ times, we see
that diam $\Delta^n_j$ is $|I^n_j|(1+O(\sum_{k=m}^n\mu_k)$. Since $\mu_k$
   decay exponentially, we are done.    \QED

\medskip\noindent
{\bf \S 3. Thurston's transformation and the shape of the complex puzzle-pieces.}
Let us consider the quadratic map $p: z\mapsto z^2-1$ and mark on ${\bf C}$
a set $A$ of three points $-1, 0,$ and $a=(1+\sqrt{5})/2$. The first two form a
cycle, while the last one is fixed.  Taking a conformal structure $\nu$ on the
thrice punctured plane $S={\bf C}\ssm A$, we can pull it back by $p$.
This induces a ``Thurston's transformation" $L $ of
 the Teichm\"{u}ller space $T_S$
of thrice punctured planes into itself
 (compare [MT] or [DH]). The main property of $L$ is
that it strictly contracts the Teichm\"{u}ller metric, and hence all
trajectories $L^n\tau$ exponentially converge to the single fixed point
$\tau_0\in T_S$ represented by the standard conformal structure.

Let us consider the involution $\rho: T_S\rightarrow T_S$
 induced by the reflection
of the conformal structure about the real line.  This involution commutes
with $L$, and so the subspace $T_S^*$ of ${\bf R}$-symmetric structures is
$L$-invariant. This subspace can be identified with the set of triples
on the real line up to affine transformations. We can normalize the
triples, say, as follows:
$\{\gamma, 0, a,\}\; \gamma<0$. To pull back such a triple, we should
take the quadratic polynomial $p_{\gamma}$
which fixes $a$ and carries 0 to $\gamma$, and take the negative preimage of
$\gamma$.

Let us rescale both intervals $I^n$ and $I^{n-1}$ to the size
$T=[-a,a]$ with $a$ as above. Let $G_n: T\rightarrow T$ be the rescaled
$g_n: I^n\rightarrow I^{n-1}$ (observe that this is a non-dynamical procedure,
compare [KP]).
Let us select the orientation in such a way that 0 is the minimum point
of $G_n$.

\proclaim Lemma 3. The maps $G_n$ converge to the polynomial $p(z)=z^2-1$
in $C^1$-norm on the compact subsets of ${\bf C}$.

\noindent {\bf Proof.}
If we pull back
the Euclidean disk $\Delta=D(I^n)$, we obtain a sequence of 
puzzle-pieces whose diameter
is commensurable with their traces on the real line (Corollary 2). 
By the Denjoy distortion argument, 
$$Dh_n^{-1}(z)=Dh_n^{-1}(0)(1+O(\sqrt{\mu_n})),\quad  z\in \Delta,$$
so that $h_n^{-1}$ in $\Delta$ is an exponentially small perturbation of
a linear map. Rescaling, we conclude that 
$G_n=H_n\circ p_{\gamma(n)}$ 
where $H_n$ are diffeomorpisms  converging exponentially
 to id in $C^1$ on compact sets, and $p_{\gamma(n)}$ are quadratic polynomials
introduced above.  

Let us now consider a sequence $\tau_n\in T^*_S$ represented by triples
$(G_n(0), 0, a)$. It was shown in  [LM] that $|G_n(0)|/a$ stays away from
0 and 1. Hence $\tau_{n+1}=L\circ Q_n (\tau_n)$
where $L$ is the Thurston transformation, while $Q_n$ is exponentially
close to id in the Teichm\"{u}ller metric. Since $L$ is strictly contracting,
$\tau_n$ must converge to its fixed point $\tau_0$. 

We conclude that $G_n(0)\to -1$, hence $p_{\gamma(n)}\to p$ and
$G_n\to p$.              \QED

Let us consider the following topology on the space ${\cal K}$
of connected compact 
subsets $K$ of ${\bf C}$. Let $\psi_K: \{z: |z|>1\}\rightarrow {\bf C}\ssm K$
be the Riemann map normalized at $\infty$ by $\psi(z)\sim qz$ with $q>0$.
Then the topology on ${\cal K}$ is induced by the compact open topology
on the space of univalent functions.

Let us now consider the complex pieces $\Delta^n$ based upon the intervals
$I^n$. Here $\Delta^n$ is the $g_n$-pull-back of $\Delta^{n-1}$. Rescaling 
of $I^n$ to $T$ leads to the corresponding rescaled  pieces $P_n$.

\proclaim Lemma 4. The pieces $P_n$ converge to the filled-in Julia set 
of $p(z)=z^2-1$.

\noindent {\bf Proof.} The piece $P_{n}$ is the $G_n$-pull-back of $P_{n-1}$.
By Lemma 1, diam$P_n$ is bounded. Hence $G_n|P_n$ 
is an exponentially small perturbation of $p$ which yields the desired.  \QED

\medskip\noindent
{\bf \S 4. Qc conjugacy on the critical sets.}
Let us consider the complex renormalizations of $f$, 
$$F_n=R^n f: V^n_0\cup V^n_1\rightarrow P^n,$$
where $V^n_i$ are the rescaled puzzle-pieces based upon the intervals
$T^n_i$. We use the same letters for the complex extensions of different maps.
In particular, let $G_n : P^n\rightarrow P^{n-1}$ is the rescaled
$g_n: \Delta^n\rightarrow \Delta^{n-1}$ (see Figure 2).  

Let us parametrize smoothly the boundary of the piece $P^0$,
$\gamma: {\bf T}\rightarrow \partial P^0$. This parametrization
can be naturally lifted to the parametrization 
$\gamma_1: {\bf T}\rightarrow\partial P^1$,
namely  $G_1\circ\gamma_1=\gamma(z^2)$,
then to the parametrization of $ \partial P^2$ etc. We refer to these
parametrizations as to the boundary markings.

Let us also consider another Fibonacci map $\tilde f$
whose data will be labeled by tilde. The {\sl Teichm\"{u}ller distance}
between two marked puzzle-pieces is the best dilatation of qc maps
between the pieces respecting the boundary marking.

\proclaim Lemma 5. The marked puzzle-pieces $P^n$ and $\tilde P^n$ 
stay bounded Teichm\"{u}ller distance apart.

\noindent{\bf Proof.} Let we have a $K$-qc map 
$H_{n-1}: P^{n-1}\rightarrow \tilde P^{n-1}$
of marked pieces respecting the positions of the critical points and the 
critical values, that is, $H_n(0)=0$ and $H_n(\gamma_n)=\tilde\gamma_n$. 
It can be lifted to the $K(1+O(\mu_n))$-qc map
$h_n:P^n\rightarrow \tilde P^n$. This map respects boundary
marking and  0-points but it does not respect $\gamma$-points.
However, it respects these points up to exponentially small error, namely
$h_n(\gamma_n$ and $\tilde\gamma_n$
are exponentially close.

  Indeed, let $q_n\in T^n_1$ be the $G_n$-preimage of 0. As the length
of $T_n$ is exponentially small, the points
$q_n$ and $\gamma_{n+1}$ are exponentially close.
 Moreover, by Lemma 4 the distance from these points  
to the boundary $\partial P^n$ is bounded from below.
By the H\"{o}lder continuity of qc maps we conclude that
$(h_n(q_n)$ and $h_n(\gamma_n)$ are also exponentially close.
As $h_n(q_n)=\tilde q_n$,  the points $h_n(\gamma_n)$
and $\tilde\gamma_n$ are exponentially close as well.

As the distance from these points to the boundary $\partial \tilde P^n$ and
from 0 is bounded from below, they are exponentially close with respect
to the Poincar\'{e} metric of $\tilde P^n$. Hence there is a diffeomorphism
$\psi: \tilde P^n\rightarrow \tilde P^n$ with exp small dilatation keeping 
$\partial \tilde P^n$ and 0 fixed, and pushing $h_n(\gamma_n)$ to
$\tilde\gamma_n$. Then $H_n=\psi\circ h_n$ is a 
$(K+$exp small)-qc map between the marked
puzzle-pieces $P_n$ and $\tilde P_n$ respecting the positions of the
critical points and the critical values.

 Proceeding in a such a way we construct  uniformly qc maps between
$P^n$ and $\tilde P^n$ on all levels
(as the exponentially small addings to dilatation sum up to a finite
value).   \QED

Let us now consider the pairs of pants $Q^n=P^n\ssm( V^n_0\cup V^n_1)$
where $V^n_0\equiv P^{n+1}$ and $V^n_1$  with naturally marked boundary. 

\proclaim Lemma 6. The pairs of pants $Q^n$ and $\tilde Q^n$ stay
bounded Teicm\"{u}ller distance apart.

\noindent{\bf Proof.} Let us consider a $K$-qc homeomorphism
$H_{n-1}: Q^{n-1}\rightarrow \tilde Q^{n-1}$
 of marked pairs of pants. It follows
from the previous lemma that we can extend these maps across $V^{n_1}_j$. 
Indeed, the previous lemma provides us with the continuation to $V^{n-1}_0$.
Moreover, it provides us with a map $P^{n-1}\rightarrow \tilde P^{n-1}$
which then can be pulled back to $V^{n-1}_1$. Let us  keep the notation 
 $H_{n-1}$ for this extension.

Let us now consider the pull-back  $W^{n-1}\subset V^{n-1}_1$ of
$V^{n-1}_0$ by $F_{n-1}$. Its boundary is also naturally marked.
By one more pull-back of $H_{n-1}$ we can reconstruct it
in such a way that it will respect this marking. 
Let us consider the annulus $A^{n-1}=P^{n-1}\ssm  W^{n-1}$
with marked boundary.

The annulus $L^n=P^n\ssm V_0^n$ double covers $A^{n-1}$ under $G_n$.
So we can pull $H^{n-1}$ back to  a $K$-qc map $H^n: L^n\rightarrow \tilde L^n$.
 Moreover, this map respects
the parametrization of $\partial  V^n_1$, and hence can be restricted
to the $K$-qc map of marked pairs of pants of level $n$.  \QED

\vskip 2in
\centerline{Figure 3}\medskip

 We are prepared to obtain the desired result of
 this section.

\proclaim Lemma 7. There is an ${\bf R}$-symmetric qc map which conjugate
 $f$ and $\tilde f$ on their critical sets.

\noindent{\bf Proof.} The critical set can be represented as
$$\omega(c)=\cap_{n=1}^{\infty} \cup Q^n_i,$$
where $Q^n_i$ are dynamically constructed disjoint pairs of pants
(see Figure 3).
They are obtained by univalent pull-backs of appropriate central
pairs of pants. As these pull-backs have bounded dilatations,
Lemma 6 implies that $Q^n_i$ stay on bounded Teichm\"{u}ller distance 
from $\tilde Q^n_i$. Gluing together all these pairs of pants, we obtain the
desired result.   \QED

\medskip\noindent{\bf \S 5. Pull-back argument.}
Sullivan's pull-back argument allows to construct a qc conjugacy between two
polynomial-like maps as long as there is a qc conjugacy on their critical
sets. In this paper we deal with asymptotically conformal maps, so that
we need the dilatation control of pull-backs.
Lemma 1 will provide us with such a control along the real line. However,
out of the real line the dilatation can grow, so that we should stop the
construction at an appropriate moment. Let us show how it works.
First we need some extra analysis on the real line.

Let $f_n:\cup I^n_j\rightarrow I^{n-1}$ be the full return map to the 
interval $I^{n-1}$.

\proclaim Lemma 8. Let $I^n\equiv J_0, J_{-1},...$ be any pull-back 
(finite or infinite) of the interval $I^n$ 
 Then $$\sum|J_{-k}|=O(\mu_n).$$

\noindent{\bf Proof.} Denote by ${\cal J}$ the union of the intervals in the
pull-back. Let us first assume that the intervals $J_{-k}$
don't intersect $I^n$. Let $K_0\equiv J_0, K_1,...$ be the piece of the
pull-back  which belongs to $I^{n-1}$, ${\cal K}={\cal J}\cap I^{n-1}$
be the union of these intervals. 
This is actually the pull-back under the map $f_n$. 
This map is expanding with bounded distortion on $I^n_j$ (actually very
strongly expanding and almost linear on $I^n_j$). Hence
$$\sum |K_{_j}|=O(\mu_n).   \eqno (6)$$
Let us now consider all intervals $L_i$ obtained by pulling $I^{n-1}$ back which
are maximal in the sense that they don't belong to another pull-back interval.
In other words, there is an  $m=m(i)$ such that $f^m L_i=I^{n-1}$ but
$f^l L_i\cap I^{n-1}=\emptyset$.  These intervals are mutually disjoint
(and cover almost everything).

Let now ${\cal K}_i={\cal J}\cap L_i$. Then $f^{m(i)}$ maps ${\cal K}_i$
with bounded distortion (actually almost linearly) onto ${\cal K}$.
Hence dens$({\cal K}_i| L_i)=O(\mu_n)$. Summing up over $i$ we get the
claim.        

Assume now that there are intervals in $I^n$ but there are no ones in
$I^{n+1}$.
 Let $J_{-l}$ be the first interval belonging to $I^n$. Then for the further
pull-back we can repeat the same argument on level $n$ instead of $n-1$
(taking into account that the Poincar\'{e} lengths of $I^{n+1}_j$ in $I^n$
are $O(\mu_n)$). 

In general case let us divide the pull-back into the pieces ${\cal J}_l$ between
the first landing at $I^l$ and the first landing at $I^{l+1}$. Let us pull
$I^l$ along the corresponding piece. This pull-back does not intersect
$I^{l+1}$ either, and according to the previous considerations its total length
is $O(\mu_l)$.
All the more this is true for the total length of ${\cal J}_l$.

Hence the total length of ${\cal J}$ is 
$O(\sum_{l\geq n}\mu_l)=O(\mu_n).$     \QED

Let us now state the complex version of the above lemma.

\proclaim Lemma 9. Let $\Omega=D(I^n), \Omega_{-1},...$ be any pull-back of
the disk $\Omega$ along the real line. Then
$$\sum {\rm diam}\Omega_{-n}=O(\mu_n).$$

\noindent {\bf Proof.} Let ${\cal W}$ denote the union of the disks in this
pull-back. As in the above argument, let us decompose it into the strings
${\cal W}_j$ in between levels $j$ and $j+1$. Let $\Omega^j$ be the first
puzzle-piece in the $j$th string. 

On the other hand, let $\Delta^j$ denote the pull-backs of $\Omega$
based upon the intervals $I^j$. Then by the  Markov property of the whole
family of pull-backs, $\Omega^j\subset \Delta^j$. 
Hence the pull-back ${\cal W}_j$
can be inscribed into the corresponding pull-back of ${\cal D}_j$ of the 
puzzle-piece $\Delta^j$.

 It follows from Lemma 1 that the sum of the diameters of pieces
in ${\cal D}_j$ is
commensurable with the total length of its trace on the real line.
By the previous lemma, the latter is $O(\mu_n)$, and we are done.  \QED

Let us now select a high level $n$ and consider the complex renormalization
$F_n : V^n_0\cup V^n_1\rightarrow P^n$. Let us re-denote all these objects
as $F: U^1_0\cup U^1_1\rightarrow U^0$. As above,  the 
corresponding objects for another Fibonacci map $\tilde f$
will be labeled with the tilde. The following statement
shows that two renormalizations of sufficiently high order are 
qc-conjugate. 

\proclaim Proposition 10.  There is a qc map $U^0\rightarrow \tilde U^0$ 
which conjugate $F$ and $\tilde F$ on the real line.

\noindent{\bf Proof.} By Lemma 7, there is  a qc map 
$h_0: U^0\rightarrow \tilde U^0$ which conjugate $F$ to $\tilde F$ on the
critical sets and on the $\partial(U^1_0\cup U^1_1)$.
Let us start to pull it back.

Let $U^n_j$ denote the family of puzzle-pieces of depth $n$ (that is, the
components of $F^{-n} U^0$) {\sl which meet the real line}. Let us assume
by induction that we have already constructed a qc map
$h_n: U^0\rightarrow \tilde U^0$ which conjugate $F$ to $\tilde F$ on their
critical sets and on $(U^1_0\cup U^1_1)\ssm{\rm int}(\cup U^n_j)$.
Then construct $h_{n+1}$ as the lift of $h_n$ to all puzzle-pieces $U^{n}_j$.

Since the puzzle-pieces $U^n_j$ shrink to points, the sequence $h_n$ has the
continuous pointwise limit $h$ which conjugate $F$ and $\tilde F$ on the
real line. Moreover,  by (2) and Lemma 9,
 the $h_n$ has uniformly bounded dilatations. Hence $h$ is qc.  \QED 

Let us re-denote $I^n$ by $J\equiv J^0$, and let $\Delta=D(J)$. 
Let us now consider the full first return map $f_1$ to $\Delta$.
 Its domain intersects the real line by the union of intervals 
$J^1_j\equiv I^{n+1}_j$. Let $\Delta^1_j$ be the pull-back of
$\Delta$ intersecting the real line by $I^{n+1}_j$,
${\cal D}^1=\cup \Delta^1_j$ (see Figure 4).

The goal of the next three lemmas is to construct a
 qc map $h: \Delta\rightarrow \tilde\Delta$
which conjugate $f_1|\partial {\cal D}$ to
 $\tilde f_1|\partial\tilde{\cal D}$
(as well as $f_1|\omega(c)$ to $\tilde f_1|\omega(\tilde c)$). This will be 
the starting data for the pull-back argument. The problem is that
the boundary $\partial{\cal D}$ is not piecewise-smooth.

Given a set $U$, denote by $U^+$  the intersection
of $U$ with the upper  half-plane.

\proclaim Lemma 11. The topological discs $\Delta^1_j$ are pairwise disjoint.
 The set $W=(\Delta\ssm {\cal D})^+$
is a quasi-disk.

\noindent{\bf Proof.} The map $f_n:\Delta^1_j\rightarrow \Delta$
has exponentially small non-linearity. Hence $\Delta^1_j$ is a minor
distorted round disk. On the other hand,  the 
intervals $J^1_i$ and $J^1_j$ are exponentially small as
compared with the gap  $G_{ij}$ in between. It follows that
the disks $\Delta^1_i$ and $\Delta_j^1$ are disjoint.

Let $\Gamma=\partial W$. It follows from the previous discussion that
this curve is
rectifiable. Take two close points $z,\zeta\in \Gamma$. Let
$\delta$ be the shortest path connecting $z$ and $\zeta$ in 
$\Gamma\cup{\bf R}$ (it is ``typically" the union of an interval of the
real line and two almost circle arcs), and $\gamma$ be the shortest arc
in $\Gamma$ connecting $z$ and $\zeta$. Then the length of $\delta$ 
is commensurable with both the length of
$\gamma$ and  the dist$(z,\zeta)$.  \QED

\vskip 2.5in
\centerline{Figure 4}\medskip

For the further discussion it is convenient to make a more special choice
of the interval $J$ (compare [GJ], [Y], [JS]). 
Namely, let $\alpha$ be the fixed point of $f$ with 
negative multiplier $\sigma\equiv f'(\alpha)$. Let ${\cal Y}^{(0)}$
be  the partition of $T$ by $\alpha$ into two intervals. Pulling this partition
back, we obtain partitions ${\cal Y}^{(n)}$ by $n$-fold preimages of $\alpha$.
Let us call the elements of this partition {\sl the puzzle-pieces of depth}
$n$. The element containing $c$ is called {\sl critical}.
We select $J=[\beta, \beta']$ 
as the critical puzzle-piece of sufficiently high depth $N$.  

Set $\tau=\log|\tilde\sigma|/\log|\sigma|$. 

Let us now start with a qc ${\bf R}$-symmetric map
$H: \Delta\rightarrow \tilde\Delta$ which carries the critical set
of $f_1$ to the critical set of $\tilde f_1$ and such that
$$ |H(z)-\tilde\beta|\asymp |z-\beta|^{\tau}. \eqno (7)$$
Moreover, let $H$ commutes with the  symmetry around $c$ induced
by $f$ and $\tilde f$. 

 Pull $H$ back to a map
$h: {\cal D}\rightarrow\tilde{\cal D}$. Since
the union $\cup J^1_j$ is dense in $J$, this map can  be continued
to a homeomorphism $h: J\rightarrow\tilde J$. Let also
 $h|\partial \Delta=H$. This defines $h$ on the topological semi-circle
$S=\partial\Delta^+$. Since $S$ and $\tilde S$ are piecewise smooth curves,
we can naturally define the notion of a quasi-symmetric map between them.

\proclaim Lemma 12. The map $h: S\rightarrow\tilde S$ is quasi-symmetric.

\noindent{\bf Proof.} Let us consider a continuation $H: T\rightarrow\tilde T$
of $H: J\rightarrow \tilde J$ which carries the puzzle-pieces of depth $N$
to the corresponding puzzle-pieces, and has the asymptotics (7) near the
boundary points of these puzzle-pieces.

Let $K$ be the expanding Cantor set of points which never land at $J$.
Each component $L$ of $T\ssm K$ (a ``gap") is a monotone pull-back of $J$
with bounded distortion. So we can pull the map $H$ back to  qs maps on
all  gaps $L$. These maps clearly glue together to a
homeomorphism $\phi: T\rightarrow \tilde T$ which respect the dynamics on the 
Cantor sets $K$ and $\tilde K$.
Moreover, if we rescale the corresponding gaps $L$ and $\tilde L$ to the
unit size then the rescaled $\phi$ near the boundary points will have 
asymptotics (7) uniformly in $L$.

Furthermore, it easily follows from the bounded distortion properties of
expanding dynamics that $\phi|K$ can be extended to a qs conjugacy $\psi$
in a neighborhood
of  $K$. This conjugacy must have the same asymtotics (7) on the rescaled
gaps (since the conjugacy near the fixed points have such asymptotics). 
It follows that $\phi$ and $\psi$ are commensurable on the gaps,
and hence $\phi$ is qs on the whole interval.

Observe now that $h: J\rightarrow\tilde J$ is the pull back of $\phi$
by the  almost quadratic maps $f|J$ and $\tilde f|\tilde J$.
Hence $h|J$ is qs and has asymptotics (7) near the boundary. Since it has
the same asymptotics on the opposite side of $\beta$, $\beta'$ on the
 arc $S\ssm J$, it is qs on $S$.  \QED

\proclaim Lemma 13. The map $h: \partial W\rightarrow \partial\tilde W$
allows a qc extension to $W\rightarrow \tilde W$.

\noindent {\bf Proof.}  Let $E$ be the exterior component of ${\bf C}\ssm S$.
By the previous lemma, there is a qc extension of $h$ from $S$ to 
$h_0:E\rightarrow \tilde E$ 
(which change the original values of $h$ below the real line).

We can now glue 
 $h: {\cal D}^+\rightarrow\tilde{\cal D}^+$ with $h_0$ to a qc map
$h_*:  {\bf C}\ssm W\rightarrow\tilde {\bf C}\ssm W$
 (since they agree on the real line).
Since $W$ is a quasi-disk (by Lemma 11), $h_*$ can be reflected to the
interior of $W$, and this is a desired extension.  \QED

\proclaim Corollary 14. There is an ${\bf R}$-symmetric qc map
 $h:\Delta\rightarrow\tilde\Delta$
which conjugates $f_1$ to $\tilde f_1$ on the critical sets and on 
the boundary of ${\cal D}$.

\noindent{\bf Proof.} Lemma 13 gives us a desired qc extension of 
the original
$h$ from ${\cal D}\cup\partial\Delta$ to $\Delta$.  \QED

     Now we are ready to prove the main result.

\proclaim Theorem I. Any two Fibonacci quasi-quadratic maps are qc
  conjugate.

\noindent {\bf Proof.} Starting with the qc map $h$ given by Corollary 14,
we can go through the pull-back argument in the same way as in Proposition 10.
This provides us with a qs conjugacy between the return maps
$f_1$ and $\tilde f_1$.  Then we can spread it around the whole interval $T$
as in the proof of  Lemma 12.  \QED

\bigskip\noindent
{\bf \S 6. Teichm\"{u}ller metric.}
Let $K_h$ denote the dilatation of a qc map $h$.
Given two Fibonacci maps $f$ and $g$ and the qs conjugacy between them,
the Teichm\"{u}ller pseudo-distance  dist$_T(f,g)$ is defined as the 
infimum of $\log K_h$ for all qc extensions of $h$.

\proclaim Theorem II. If {\rm dist}$_T(f,g)=0$ then $f$ and $g$ are smoothly
  conjugate.

\noindent{\bf Proof.} Our first step is the same as Sullivan's [S1]:
If dist$_T(f,g)=0$ then the multipliers of the 
corresponding periodic orbits of the maps are equal. However, as we
don't have yet a proper thermodynamical formalism for unimodal maps,
we will proceed by a concrete geometric analysis.

The next observation is that the parameter $a$ in (1) must be the same 
for $f$ and $g$. Indeed, it can be explicitly expressed via
 the multipliers
of the fixed points of the return maps $g_n: I^n\rightarrow I^{n-1}$ 
(since the $g_n$ are asymptotically quadratic). By [LM] this already yields
the smoothness of the conjugacy on the critical sets.

Let us now take a point $x\in I^n\ssm I^{n-1}$ and push it forward by iterates 
of $g_n$ till the first moment it lands in $I^n$ (if any), then apply the
iterates of $g_{n+1}$ till the first moment it lands in $I^{n+1}$, etc.
This provides us with a nested sequence of intervals around $x$ whose
lengths can be expressed (up to a bounded error) through the scaling
factors and the multipliers of appropriate periodic points (by shadowing).
This implies that $h$ is Lipschitz continuous.
Moreover, when we approach the critical point, then the errors in the
above argument exponentially decrease. Hence $h$ is smooth at the 
critical point.

Given now any pair of intervals $I\supset J$, let us show that
$$\left|{|h J|\over |J|} : {|hI|\over |I|}-1\right|=O(|I|).  \eqno (8) $$
This is enough to prove locally at any point $a$.
By the previous considerations, this is true
at the critical point. Since the critical set $\omega(c)$ is minimal, 
this is also true for any $a\in \omega(c)$.

Let now $a\not\in\omega(c)$, and $I$ be a tiny interval around $I$.
Remark that almost all ponts $x\in I$ eventually  return back to $I$.
Let us take the pull-back of $I$ corresponding to this return.
This provides us with the covering of almost all of $I$ by intervals $L_k$.
The distortion of the return map  $g$ is $O(|I|
)$ on the all $L_k's$.
Let $\sigma_k$ be the multiplier of the $g$-fixed point in $L_k$.
Then we conclude that
$$\left|{|I|\over |L_k|}:\sigma_k-1\right|=O(|I|), \eqno(9)  $$
and the analogous estimate holds for the second map. Since the corresponding
multipliers of these maps are equal, we obtain (8) with $J=L_k$.
Repeating now this procedure for returns of higher order, we obtain
an arbitrarily fine covering of almost the whole of $I$ by intervals
for which (8) hold. This implies (8) for any $J\subset I$.

Let  $\epsilon_n=1/2^n$, and let us consider the sequence of functions
$$\rho_n(x)={h(x+\epsilon_n)-h(x-\epsilon_n)\over 2\epsilon_n}. $$
According to (8) and Lipschitz continuity
$$|\rho_n(x)-\rho_{n+1}(x)|=O(\epsilon_n)  \eqno (10) $$
uniformly in $x$. Hence the $\rho_n$ uniformly converge to the derivative 
of $h$.

\bigskip
\centerline{\bf References.}\smallskip

{\item [DH]} A.Douady \& J.H.Hubbard. A proof of Thurston's topological
   characterization of rational functions. Preprint Institut
   Mittag-Leffler, 1986.\smallskip

{\item [G]} F.P.Gardiner. Lacunary series as quadratic differentials,
   Preprint 1992.   

{\item[GJ]} J. Guckenheimer \& S. Johnson. Distortion of S-unimodal maps.
   Annals of Math., {\bf 132} (1990), 71-130. \smallskip

{\item [H]} J.H. Hubbard.  Local connectivity of Julia sets and bifurcation
     loci: three theorems of J.-C. Yoccoz. In:
    ``Topological Methods in Modern Mathematics, A Symposium in
     Honor of John Milnor's 60th Birthday", Publish or Perish, 1993.\smallskip

{\item [HK]} F.Hofbauer and G.Keller. Some remarks on recent results 
   about S-unimodal maps. Ann. Institut Henri Poincar\'{e}, {\bf 53} (1990),
  413-425.\smallskip

{\item [Ji]} Y.Jiang. Generalized Ulam-von Neumann transformations.
    Thesis, 1990.\smallskip

{\item [JS]} M.Jakobson and G.Swiatek. Quasisymmetric conjugacies
   between unimodal maps.  Preprint IMS at Stony Brook, \#1991/16.
   \smallskip

{\item [KN]} G.Keller \& T.Nowicki. Fibonacci maps revisited.
             Preprint, 1992. \smallskip

{\item [KP]} J.Ketoja \& O.Piirila. On the abnormality of the period
doubling bifurcation. Phys. Letters A, {\bf 138} (1989), 488 - 492. \smallskip 

{\item [K]} J. Kahn. Holomorphic Removability of Julia Sets. Manuscript in
  preparation.\smallskip

{\item [LM]} M.Lyubich \& J.Milnor. The unimodal Fibonacci map. 
        Preprint IMS  at Stony Brook \#1991/15. To appear in the Journal
        of AMS. \smallskip

{\item [L]} M.Lyubich. Combinatorics, geometry and attractors of
     quasi-quadratic maps.    
     Preprint IMS  at Stony Brook \#1992/18.\smallskip

{\item [MT]} J.Milnor \& W.Thurston. On iterated maps of the interval,
Preprint of 1977 and
pp. 465-563 of Dynamical Systems, Proc. U. Md., 1986-87,
ed. J. Alexander, Lect. Notes Math., 1342, Springer 1988.\smallskip

{\item [S1]} D.Sullivan. Quasiconformal homeomorphisms in dynamics, topology
    and geometry.   Proceedings of the ICM, Berkeley, 1986, {\bf 2},
    1216.

{\item [S2]} D.Sullivan. Bounds, quadratic differentials, and renormalization
conjectures,  1990. To appear in AMS Centennial Publications.
{\bf 2}: Mathematics into Twenty-first Century. \smallskip

\end